\documentclass[11pt]{article}
\usepackage{graphicx}
\usepackage{amsmath}
\usepackage{amssymb}
\usepackage{amsthm}
\usepackage[titletoc]{appendix}
\usepackage{url}
\newcommand\al\alpha
\newcommand\be\beta
\newcommand\de\delta
\newcommand\ep\varepsilon
\newcommand\tha\theta
\newcommand\ka\kappa
\newcommand\la\lambda
\newcommand\om\omega

\newcommand\iy\infty
\newcommand\pa\partial

\newcommand{\hyp}[5]{\,\mbox{}_{#1}F_{#2}\!\left(\genfrac{}{}{0pt}{}{#3}{#4};#5\right)}

\numberwithin{equation}{section}

\newtheorem{theorem}{Theorem}

\newtheorem{Remark}[theorem]{Remark}

\begin{document}

\title{The Catalan number of the second kind and a related integral.}
\author{Enno Diekema \footnote{email adress: e.diekema@gmail.com}}
\maketitle

\begin{abstract}
\noindent
Li et al. give an integral formula for the Catalan number of the second kind. They show that this integral can be written as a summation with double factorials. In this paper the integral is reduced to a product of the Catalan number and a hypergeometric function. This hypergeometric function can be written as an associated Legendre polynomial of the first kind and as a Jacobi polynomial which is however not orthogonal. In the last section a related integral is calculated.   
\end{abstract}

\section{Introduction}
\setlength{\parindent}{0cm}
The Catalan numbers are a sequence of positive integers that appear in many counting problems in combinatorics. They are named after the mathematician E.C. Catalan (1814-–1894). They have a closed-form formula in terms of factorials. They also satisfy a fundamental recurrence relation. Explicit formula for the Catalan numbers $C_n$ are 
\begin{equation}
C_n=\dfrac{(2n)!}{n!(n+1)!}=\dfrac{1}{n+1}\binom{2n}{n}=\dfrac{4^n\ \Gamma(n+1/2)}{\sqrt{\pi}\ \Gamma(n+2)}=
\hyp21{-n,1-n}{2}{1}.
\label{0}
\end{equation}
There are a lot more. The first Catalan numbers for $0 \leq n \leq 8$ are:
\[
C_0..._8=1,\, 1,\, 2,\, 5,\, 14,\, 42,\, 132,\, 429.
\]

A simple recurrence equation is
\[
C_{n+1}=\dfrac{2(2n+1)}{n+2}C_n \quad \text{with}\quad C_0=1.
\]
Among the different generation functions we take
\[
G(x)=\dfrac{2}{1+\sqrt{1-4x}}=\sum_{n=0}^{\infty}C_n\ x^n.
\] 

There are many generalizations of the Catalan numbers. Following the previous formula it is a natural generalization to take the generation function \cite{1}, \cite{2}
\[
G_{a,b}(x)=\dfrac{1}{a+\sqrt{b-x}}=\sum_{n=0}^\infty\mathcal{C}_n(a,b)x^n.
\]
for $a \geq 0$ and $b>0$. We get $\mathcal{C}_n(1/2,1/4)=C_n$. The authors in \cite{1} call this generalization the {\em Catalan numbers of the second kind}. They derive a formula with a summation with double factorials for these numbers \cite[2.1]{1}
\[
\mathcal{C}_n(a,b)=\dfrac{1}{(2n)!!\ b^{n+1/2}}\sum_{k=0}^n\binom{2n-k-1}{2(n-k)}
\dfrac{k!\ [2(n-k)-1]!!}{(1+a/\sqrt{b})^{k+1}}.
\]

They also give an integral formula for these numbers with $a,b>0$ \cite[4.1]{3}
\begin{equation}
\mathcal{C}_n(a,b)=\dfrac{1}{\pi}
\int_0^\infty\dfrac{\sqrt{t}}{(a^2+t)}\dfrac{1}{(b+t)^{n+1}}dt, \qquad n \geq 0.
\label{0a}
\end{equation}
The aim of the first part of this paper is to determine a closed form of this integral. The main result of the first part of this paper is with $a \geq 0 \ \ b,n>0$
\begin{align}
\mathcal{C}_n(a,b)
&=C_n\dfrac{\pi}{(2\sqrt{b})^n}\dfrac{1}{\big(a+\sqrt{b}\big)^{n+1}}
\hyp21{1-n,n}{n+2}{\dfrac{\sqrt{b}-a}{2\sqrt{b}}} \nonumber \\
&=C_n\left(\dfrac{a}{2\sqrt{b}}\right)^n
\dfrac{\Gamma(n+2)}{\big(\sqrt{b}-a\big)^{2n+1}}P^{-n-1}_{n-1}\left(\dfrac{a}{\sqrt{b}}\right) \nonumber \\
&=\dfrac{\pi}{n\big(2\sqrt{b}\big)^n}\dfrac{1}{\big(a+\sqrt{b}\big)^{n+1}}
P_{n-1}^{(n+1,-n-1)}\left(\dfrac{a}{\sqrt{b}}\right)
\label{0b}
\end{align}
where $_2F_1$ is a hypergeometric function. $P^{-n-1}_{n-1}(x)$ is an associated Legendre function and $P_{n-1}^{(n+1,-n-1)}(x)$ is a Jacobi polynomial which is however not orthogonal.

\

In \cite{2} the authors define a Catalan functional which is given by
\begin{equation}
\mathfrak{C}_n(a,b;p)=\dfrac{\sin(p\pi)}{\pi}
\int_0^\infty\dfrac{t^p}{a^2+2a\cos(p\pi)\, t^p+t^{2p}}\dfrac{1}{(b+t)^{n+1}}dt
\label{B1}
\end{equation}
with conditions:\ $a \geq 0,\ b>0,\ 0<p<1,\ \ n \geq 0$.

They also give a double summation as the outcome of the integral \cite[(5.1)]{2}
\begin{multline}
\int_0^\infty\dfrac{t^p}{a^2+2a\cos(p\pi)\, t^p+t^{2p}}\dfrac{1}{(b+t)^{n+1}}dt= \\
=\dfrac{\pi}{(a+b^p)\sin(p\pi)}\dfrac{(-1)^n}{\Gamma(n+1)b^n}\sum_{k=0}^n
\dfrac{1}{\big(1+a/b^p\big)^k}\sum_{m=0}^k(-1)^m\binom{k}{m}\langle p\, m\rangle_n
\label{B1a}
\end{multline}
with $\langle a\, \rangle_n$ the falling factorial. Because there is often a difference between the meaning of the rising and falling factorials in combinatoric theory and the special function theory we give here the definitions from the theory of special functions. The falling factorial is defined by
\[
\langle a\, \rangle_n=
\begin{cases}
a(a-1)\dots\big(a-(n-1)\big) \qquad n \geq 1 \\
1, \qquad\qquad\qquad \qquad\qquad\quad\ \ n=0.
\end{cases}
\]
Then we can write
\[
\langle a\, \rangle_n=\dfrac{a(a-1)\dots\big(a-(n-1)\big)\Gamma\big(a-(n-1)\big)}{\Gamma\big(a-(n-1)\big)}=\dfrac{\Gamma(a+1)}{\Gamma(a+1-n)}.
\]
However in the theory of the special functions one more commonly uses the rising factorial (also called Pochhammer symbol). This is defined as
\[
(a)_n=a(a+1)\dots(a+n-1).
\]
This can be written as
\[
(a)_n=\dfrac{a(a+1)\dots(a+n-1)\Gamma(a)}{\Gamma(a)}=\dfrac{\Gamma(a+n)}{\Gamma(a)}.
\]
For the connection between the falling factorial and the rising factorial there is
\[
\langle a\, \rangle_n=(-1)^n(-a)_n.
\]

The aim of the second part of this paper is to develop the integral with the presumption that we look for a solution with a single summation. The main result of the second part of this paper is 
\[
\mathfrak{C}_n(a,b;p)=\pm\dfrac{1}{a\, b^n(n+1)}
\sum_{k=0}^\infty(\mp p\, k)_n\left(-\dfrac{b^p}{a}\right)^{\pm k}
\qquad \left|\dfrac{b^p}{a}\right|^{\pm 1}<1.
\]
Much more representations are given. For $b^p=a$ one has to take the limit. In that case care should be taken for $n=0$.\

\section{Conversion of the Catalan numbers of the second kind into a bounded hypergeometric function}

To work out the integral \eqref{0a} we use \cite[2.1.3 (12)]{8}.
\[
\int_0^\infty s^{\beta-1}(1+s)^{\gamma-\beta-1}(1+s\, z)^{-\alpha}ds=
\dfrac{\Gamma(\beta)\Gamma(\alpha+1-\gamma)}{\Gamma(\alpha+\beta-\gamma+1)}
\hyp21{\alpha,\beta}{\alpha+\beta-\gamma+1}{1-z}
\]

Application to \eqref{0a} gives
\[
\mathcal{C}_n(a,b)=\dfrac{\sqrt{\pi}\Gamma(n+1/2)}{2\Gamma(n+2)}\dfrac{b^{1/2-n}}{a^2}
\hyp21{1,3/2}{n+2}{1-\dfrac{b}{a^2}}.
\]

Using the definition of the Catalan number \eqref{0} gives
\begin{equation}
\mathcal{C}_n(a,b)=C_n\dfrac{\pi}{2^{2n+1}}\dfrac{b^{1/2-n}}{a^2}
\hyp21{1,3/2}{n+2}{1-\dfrac{b}{a^2}}.
\label{6}
\end{equation}
The Catalan number of the second kind can be written as a product of the Catalan number and a hypergeometric function. The hypergeometric function gives a polynomial. Because it is not clear that this hypergeometric function is bounded, it is usefull to convert this hypergeometric function into a hypergeometric function with a negative integer upper parameter. We use a quadratic transformation of the hypergeometric function \cite[15.3.19]{4}
\[
\hyp21{1,3/2}{n+2}{1-y^2}
=\dfrac{4}{(1+y)^2}\hyp21{1-n,2}{n+2}{\dfrac{1-y}{1+y}}
=\dfrac{1}{y^2}\left(\dfrac{2y}{1+y}\right)^{n+1}
\hyp21{1-n,n}{n+2}{\dfrac{y-1}{2y}}.
\]
Application to \eqref{6} gives at last
\begin{equation}
\mathcal{C}_n(a,b)=C_n\dfrac{\pi}{(2\sqrt{b})^n}\dfrac{1}{\big(a+\sqrt{b}\big)^{n+1}}
\hyp21{1-n,n}{n+2}{\dfrac{\sqrt{b}-a}{2\sqrt{b}}}   \qquad  a \geq 0 \ \ b,n>0. 
\label{7}
\end{equation}
The hypergeometric function can converted into an associated Legendre function of the first kind. We use first a Gauss transformation on \eqref{6}
\[
\mathcal{C}_n(a,b)=C_n\dfrac{\pi}{(4b)^{n+1/2}}\hyp21{1,n+1/2}{n+2}{1-\dfrac{a^2}{b}}.
\]
Then we use \cite[7.3.1. (41)]{3}
\[
\hyp21{a,b}{a+b+1/2}{x}=2^{a+b-1/2}\Gamma\left(a+b+\dfrac{1}{2}\right)
x^{(1-2a-2b)/4}P^{1/2-a-b}_{a-b-1/2}\big(\sqrt{1-x}\big) \qquad 0<x<1.
\]
Application with \cite[8.2.1]{1} gives
\[
\mathcal{C}_n(a,b)
=C_n\dfrac{\pi}{\big(2\sqrt{b}\big)^n}
\dfrac{\Gamma(n+2)}{\big(\sqrt{b-a^2}\big)^{n+1}}P^{-n-1}_{n-1}\left(\dfrac{a}{\sqrt{b}}\right)   \qquad  a \geq 0 \ \ b>0. 
\]

\

The hypergeometric function can converted into a Jacobi polynomial. Using \cite[22.5.42]{1} gives
\[
\mathcal{C}_n(a,b)=
\dfrac{\pi}{n\big(2\sqrt{b}\big)^n}\dfrac{1}{\big(a+\sqrt{b}\big)^{n+1}}
P_{n-1}^{(n+1,-n-1)}\left(\dfrac{a}{\sqrt{b}}\right)
  \qquad  a \geq 0 \ \ b,n>0.
\]

For the first six Catalan numbers of the second kind we find
\begin{align*}
&\mathcal{C}_0(a,b)=\pi\dfrac{1}{a+\sqrt{b}},
\qquad\qquad\qquad\qquad\qquad\qquad\qquad\qquad\qquad\qquad\qquad\qquad\qquad\qquad \\
&\mathcal{C}_1(a,b)=\pi\dfrac{1}{2(a+\sqrt{b})^2\ b^{1/2}}, \\
&\mathcal{C}_2(a,b)=\pi\dfrac{a+3\sqrt{b}}{8(a+\sqrt{b})^3\ b^{3/2}}, \\
&\mathcal{C}_3(a,b)=\pi\dfrac{a^2+4a\sqrt{b}+5b}{16(a+\sqrt{b})^4\ b^{5/2}}, \\
&\mathcal{C}_4(a,b)=\pi\dfrac{5a^3+25a^2\sqrt{b}+47ab+35b^{3/2}}{128(a+\sqrt{b})^5\ ^{7/2}}, \\
&\mathcal{C}_5(a,b)=\pi\dfrac{7a^4+42a^3\sqrt{b}+102a^2b+122ab^{3/2}+63b^2}{256(a+\sqrt{b})^6\ b^{9/2}}.
\end{align*}

\section{The Catalan functional}
In this section we derive a simple formula for the Catalan functional. Combining \eqref{B1} and \eqref{B1a} gives a formula for the Catalan functional
\begin{multline*}
\mathfrak{C}_n(a,b;p)=\dfrac{\sin(p\pi)}{\pi}\int_0^\infty\dfrac{t^p}{a^2+2a\cos(p\pi)\, t^p+t^{2p}}\dfrac{1}{(b+t)^{n+1}}dt= \\
=\dfrac{1}{(a+b^p)}\dfrac{(-1)^n}{\Gamma(n+1)b^n}\sum_{k=0}^n
\dfrac{1}{\big(1+a/b^p\big)^k}\sum_{m=0}^k(-1)^m\binom{k}{m}\langle p\, m\rangle_n.
\end{multline*}

Using the connection between the falling and the rising factorial we get
\begin{multline}
\mathfrak{C}_n(a,b;p)=\dfrac{\sin(p\pi)}{\pi}\int_0^\infty\dfrac{t^p}{a^2+2a\cos(p\pi)\, t^p+t^{2p}}\dfrac{1}{(b+t)^{n+1}}dt= \\
=\dfrac{1}{(a+b^p)}\dfrac{1}{\Gamma(n+1)b^n}\sum_{k=0}^n
\dfrac{1}{\big(1+a/b^p\big)^k}\sum_{m=0}^k(-1)^m\binom{k}{m}(-p\, m)_n.
\label{1a}
\end{multline}

For $n=0$ there is a special case. There remains
\begin{equation}
\mathfrak{C}_0(a,b;p)=\dfrac{1}{a+b^p}.
\label{1b}
\end{equation}

\

The integral can be calculated with partial integration. Therefore we simplify the denominator of the first fraction
\[
a^2+2a\cos(p\pi)\, t^p+t^{2p}={\big(t^p+a\, \cos(p\pi)\big)^2+a^2 \sin^2(p\pi)}=
(t^p+c)^2+d^2
\]
with
\[
c=a\, \cos(p\pi)\qquad \text{and}\qquad d=a\, \sin(p\pi)
\]
and this results in
\[
c+i\, d=a\,e^{ip\pi}\qquad \text{and}\qquad c-i\, d=a\, e^{-ip\pi}.
\]
Setting
\[
\dfrac{d}{dt}u(t)=\dfrac{t^p}{(t^p+c)^2+d^2}\qquad \text{and}\qquad
v(t)=\dfrac{1}{(b+t)^{n+1}}
\]
gives
\begin{align}
u(t)&=\int\dfrac{t^p}{(t^p+c)^2+d^2}dt= \nonumber \\
&=-i\dfrac{e^{ip\pi}}{2d\, a(p+1)}t^{p+1}
\hyp21{1,1+1/p}{2+1/p}{-\dfrac{e^{ip\pi}}{a}t^p}
+i\dfrac{e^{-ip\pi}}{2d\, a(p+1)}t^{p+1}
\hyp21{1,1+1/p}{2+1/p}{-\dfrac{e^{-ip\pi}}{a}t^p}
\label{B2}
\end{align}
and 
\[
\dfrac{d}{dt}v(t)=-\dfrac{n+1}{(b+t)^{n+2}}.
\] 
The function $u(t)$ can be simplified. Writing the hypergeometric functions as summations gives
\begin{multline*}
u(t)=-i\dfrac{1}{2d\, a(p+1)}
\sum_{k=0}^\infty \dfrac{(1+1/p)_k}{(2+1/p)_k}\left(-\dfrac{1}{a}\right)^k e^{ip(k+1)\pi}t^{p\,k+p+1}+ \\
+i\dfrac{1}{2d\, a(p+1)}
\sum_{k=0}^\infty \dfrac{(1+1/p)_k}{(2+1/p)_k}\left(-\dfrac{1}{a}\right)^k e^{-ip(k+1)\pi}t^{p\,k+p+1}.
\end{multline*}
Combining the terms gives
\begin{align*}
u(t)
&=\dfrac{1}{d\, a(p+1)}
\sum_{k=0}^\infty \dfrac{(1+1/p)_k}{(2+1/p)_k}\left(-\dfrac{1}{a}\right)^k t^{p\, k+p+1}
\left(\dfrac{e^{ip(k+1)\pi}-e^{-ip(k+1)\pi}}{2i}\right) \\
&=\dfrac{1}{d\, a(p+1)}
\sum_{k=0}^\infty \dfrac{(1+1/p)_k}{(2+1/p)_k}\left(-\dfrac{1}{a}\right)^k t^{p\, k+p+1}
\sin\big(p(k+1)\pi\big) \\
&=-\dfrac{\pi}{d\, a}
\sum_{k=0}^\infty \dfrac{1}{\Gamma(2+p+p\, k)\Gamma(-p-p\, k)}
\left(-\dfrac{1}{a}\right)^k t^{p\, k+p+1} .
\end{align*}
Partial integration gives
\begin{align*}
\int_0^\infty\dfrac{d}{dt}u(t)v(t)dt
&=\Big[u(t)v(t)\Big]_0^\infty-\int_0^\infty u(t)\dfrac{d}{dt}v(t)dt \\
&=\left[-\dfrac{\pi}{d\, a}
\sum_{k=0}^\infty \dfrac{1}{\Gamma(2+p+p\, k)\Gamma(-p-p\, k)}
\left(-\dfrac{1}{a}\right)^k\dfrac{t^{p\, k+p+1}}{(b+t)^{n+1}}\right]_0^\infty+ \\
&+\dfrac{\pi(n+1)}{d\, a}
\sum_{k=0}^\infty\dfrac{1}{\Gamma(2+p+p\, k)\Gamma(-p-p\, k)}
\left(-\dfrac{1}{a}\right)^k \int_0^\infty\dfrac{t^{p\, k+p+1}}{(b+t)^{n+2}}dt.
\end{align*}
Here we have interchanged the integral and the summation which is allowed because the summation is convergent. The terms within the square brackets becomes zero. For the integral we get
\[
\int_0^\infty\dfrac{t^{p\, k+k+1}}{(b+t)^{n+2}}dt=
\Gamma(2+p+p\, k)\Gamma(n-p-p\, k)\dfrac{b^{-n+p+p\, k}}{\Gamma(n+2)}.
\]
Substitution gives
\[
\mathfrak{C}_n(a,b;p)=\dfrac{\sin(p\pi)}{d\, a\, b^{n-p}(n+1)}
\sum_{k=0}^\infty\dfrac{\Gamma(n-p-p\, k)}{\Gamma(-p-p\, k)}\left(-\dfrac{b^p}{a}\right)^k.
\]
Using $d=a\sin(p\pi)$ and writing the quotient of the Gamma functions as a Pochhammer symbol gives
\[
\mathfrak{C}_n(a,b;p)=\dfrac{1}{a^2\, b^{n-p}(n+1)}
\sum_{k=0}^\infty(-p-p\, k)_n\left(-\dfrac{b^p}{a}\right)^k.
\]
Shifting the summation variable $k+1$ to $k$ results at last 
\begin{equation}
\mathfrak{C}_n(a,b;p)=\dfrac{1}{a\, b^n(n+1)}
\sum_{k=0}^\infty(-p\, k)_n\left(-\dfrac{b^p}{a}\right)^k.
\label{B3}
\end{equation}
The convergence condition for the summation is\ \  $\left|\dfrac{b^p}{a}\right|<1$.

\

\

If $\left|\dfrac{b^p}{a}\right|>1$ we can transform the hypergeometric functions in \eqref{B2} with \cite[15.3.7]{4}. Application gives
\begin{align*}
t^{p+1}e^{ip\pi}\hyp21{1,1+1/p}{2+1/p}{-\dfrac{e^{ip\pi}}{a}t^p}
&=a(p+1)\, t\, \hyp21{1,-1/p}{1-1/p}{-a\,e^{-ip\pi}t^{-p}}+ \\
&+a(p+1)\Gamma(1+1/p)\Gamma(1-1/p)\, a^{1/p}e^{-i\pi}
\end{align*}
\begin{align*}
t^{p+1}e^{-ip\pi}\hyp21{1,1+1/p}{2+1/p}{-\dfrac{e^{-ip\pi}}{a}t^p}
&=a(p+1)\, t\, \hyp21{1,-1/p}{1-1/p}{-a\,e^{ip\pi}t^{-p}}+ \\
&+a(p+1)\Gamma(1+1/p)\Gamma(1-1/p)\, a^{1/p}e^{i\pi}.
\end{align*}
Substitution in \eqref{B2} gives
\begin{align*}
u(t)&=\int\dfrac{t^p}{(t^p+c)^2+d^2}dt= \\
&=\dfrac{i}{2d}
t\, \hyp21{1,-1/p}{1-1/p}{-a\,e^{ip\pi}t^{-p}}
-\dfrac{i}{2d}
t\, \hyp21{1,-1/p}{1-1/p}{-a\,e^{-ip\pi}t^{-p}} + \\
&+i\dfrac{\Gamma(1+1/p)}{2d}\Gamma(1-1/p)\, a^{1/p}e^{i\pi}
-i\dfrac{\Gamma(1+1/p)}{2d}\Gamma(1-1/p)\, a^{1/p}e^{-i\pi} \\
&=\dfrac{i}{2d}\, t\, 
\left(\hyp21{1,-1/p}{1-1/p}{-a\,e^{ip\pi}t^{-p}}
-\hyp21{1,-1/p}{1-1/p}{-a\,e^{-ip\pi}t^{-p}}\right).
\end{align*}
Writing the hypergeometric functions as summations gives
\[
u(t)=-\dfrac{1}{d} \sum_{k=0}^\infty \dfrac{(-1/p)_k}{(1-1/p)_k}(-a)^k \left(\dfrac{e^{ip\, k\pi}-e^{-ip\, k\pi}}{2i}\right)t^{1-p\, k}.
\]
Simplifying gives
\[
u(t)=-\dfrac{\pi}{d} 
\sum_{k=0}^\infty \dfrac{1}{\Gamma(p\, k)\Gamma(2-p\, k)}(-a)^k t^{1-p\, k}.
\]
Partial integration gives
\begin{align*}
\int_0^\infty\dfrac{d}{dt}u(t)v(t)dt
&=\Big[u(t)v(t)\Big]_0^\infty-\int_0^\infty u(t)\dfrac{d}{dt}v(t)dt \\
&=-\dfrac{\pi}{d} 
\sum_{k=0}^\infty \dfrac{1}{\Gamma(p\, k)\Gamma(2-p\, k)}(-a)^k \left[\dfrac{t^{1-p\, k}}{(b+t)^{n+1}}\right]_0^\infty- \\
&-\dfrac{\pi(n+1)}{d} 
\sum_{k=0}^\infty \dfrac{1}{\Gamma(p\, k)\Gamma(2-p\, k)}(-a)^k 
\int_0^\infty \dfrac{t^{1-p\, k}}{(b+t)^{n+2}}dt. 
\end{align*}
We have interchanged the integral and the summation which is allowed because the summation is convergent. The terms with the square brackets becomes zero. For the integral we get
\[
\int_0^\infty \dfrac{t^{1-p\, k}}{(b+t)^{n+2}}dt=
\dfrac{1}{b^{n+p\, k}}\dfrac{\Gamma(2-p\, k)\Gamma(n+p\, k)}{\Gamma(n+2)}.
\]
Substitution gives
\[
\int_0^\infty\dfrac{d}{dt}u(t)v(t)dt=
-\dfrac{\pi}{d\, b^n\Gamma(n+1)} 
\sum_{k=0}^\infty (p\, k)_n\left(-\dfrac{a}{b^p}\right)^k 
\qquad \left|\dfrac{a}{b^p}\right|<1.
\]
Substitution this integral with $d=a\, \sin(p\pi)$ in \eqref{B1} gives
\[
\mathfrak{C}_n(a,b;p)=-
\dfrac{1}{a\, b^n\Gamma(n+1)} 
\sum_{k=0}^\infty (p\, k)_n\left(-\dfrac{a}{b^p}\right)^k 
\qquad \left|\dfrac{a}{b^p}\right|<1.
\]
Combining with \eqref{B3} gives
\begin{equation}
\mathfrak{C}_n(a,b;p)=\pm\dfrac{1}{a\, b^n(n+1)}
\sum_{k=0}^\infty(\mp p\, k)_n\left(-\dfrac{b^p}{a}\right)^{\pm k}
\qquad \left|\dfrac{b^p}{a}\right|^{\pm 1}<1.
\label{4.3}
\end{equation}
For $n>0$ and $\left|\dfrac{b^p}{a}\right|=1$ we had to take the limit. For the case $n=0$ and $b^p=a$ we get from \eqref{1b}
\[
\mathfrak{C}_0(b^p,b;p)=\dfrac{1}{2a}=\dfrac{1}{2b^p}.
\]

There is an important difference between the formulas \eqref{B1a} and \eqref{4.3}. In \eqref{B1a} there are two summations which are both bounded. Formula \eqref{4.3} has a single summation but the upper bound is infinity. So formula \eqref{4.3} is less suitable for numerical computing the approximated value of the integral. 

\ 

To develop other possibilities for \eqref{4.3} we set $\dfrac{b^p}{a}=y$ and define the function
\begin{equation}
Q(n,y,p)=\sum_{k=0}^\infty(-p\, k)_n(-y)^k \qquad 0 \leq y \leq 1 \qquad 0<p<1.
\label{4.4}
\end{equation}
So for \eqref{1a} we get
\[
\mathfrak{C}_n(a,b;p)=\dfrac{y}{b^{n+p}(n+1)}Q(n,y,p) \qquad 0 \leq y\leq 1.
\]
For $y=0$ we get:
\begin{align*}
Q(n,0,p)&=1 \qquad n=0 \\
&=0 \qquad n>0.
\end{align*}

In the next section we develop a number of properties of the function $Q(n,y,p)$.

\section{The function $Q(n,y,p)$}
The function $Q(n,y,p)$ is a summation of terms with fractions. The nominator of these fractions is a polynomial in $p$ of degree $n$. The denominator of these fractions are powers of $(y+1)^n$. For $n=0...4$ we have the following table

\begin{align*}
&Q(0,y,p)=\dfrac{1}{(y+1)}, \\
&Q(1,y,p)=-\dfrac{p}{(y+1)^2}+\dfrac{p}{(y+1)}, \\
&Q(2,y,p)=\dfrac{2p^2}{(y+1)^3}-\dfrac{p+3p^2}{(y+1)^2}+\dfrac{p+p^2}{(y+1)}, \\ 
&Q(3,y,p)=-\dfrac{6p^3}{(y+1)^4}+\dfrac{6(p^2+2p^3)}{(y+1)^3}-
\dfrac{2p+9p^2+7p^3}{(y+1)^2}+\dfrac{2p+3p^2+p^3}{(y+1)}, \\
&Q(4,y,p)=\dfrac{24p^4}{(y+1)^5}-\dfrac{12(3p^3+5p^4)}{(y+1)^4}+
\dfrac{2(11p^2+36p^3+25p^4)}{(y+1)^3}-\dfrac{3(2p+11p^2+14p^3+5p^4)}{(y+1)^2}+ \\
&\quad\quad\quad\ \ \ 
+\dfrac{(6p+11p^2+6p^3+p^4)}{(y+1)}.
\end{align*}
It is clear that the nominators of the terms in $Q(n,y,p)$ are polynomials of degree $n$. So $Q(n,y,p)$ can be written as a double summation. But the form is different from that of \eqref{B1a}. One recognizes in the first terms the factorial of $n$. In the last term there appears a polynomial with Stirling numbers of the first kind.

When adding the terms we can rearrange them. Setting $z=y+1$ we get

\begin{align*}
&Q(0,z,p)=\dfrac{1}{z}, \\
&Q(1,z,p)\left[p\dfrac{z-1}{z^2}\right]^{-1}=1, \\
&Q(2,z,p)\left[p\dfrac{z-1}{z^3}\right]^{-1}=z+p(z-2), \\
&Q(3,z,p)\left[p\dfrac{z-1}{z^4}\right]^{-1}=2z^2+3\,p\,z(z-2)+p^2(z^2-6z+6), \\
&Q(4,z,p)\left[p\dfrac{z-1}{z^5}\right]^{-1}=6z^3+11\,p\,z^2(z-2)+6\,p^2z(z^2-6z+6)+p^3(z^3-14z+36z-24), \\
&Q(5,z,p)\left[p\dfrac{z-1}{z^6}\right]^{-1}=24z^4+50\,p\,z^3(z-2)+35\,p^2z^2(z^2-6z+6)+10\,p^3z(z^3-14z+36z-24)+ \\
&\qquad\qquad\qquad\qquad\quad\  +p^4(z^4-30z^3+150z^2-240z+120).
\end{align*}
The coefficients of the polynomials in brackets are similar to the coefficients of the geometric polynomials. The geometric polynomials $\omega_n(x)$ and their inverses are defined as \cite{5},\cite[Problem 12,\ pp 310]{6}
\[
\omega_n(x)=\sum_{k=0}^n S(n,k)k!x^k\qquad \text{and} \qquad x^n=\dfrac{1}{n!}\sum_{k=0}^n(-1)^{n-k}s(n,k)\omega_k(x)
\]
with $s(k,m)$ the Stirling numbers of the first kind and $S(k,m)$ the Stirling numbers of the second kind. For $n=0,1,\dots,5$ we get
\begin{align*}
&\omega_0(x)=1,\qquad \omega_1(x)=x,\qquad \omega_2(x)=x(2x+1),\qquad \omega_3(x)=x(6x^2+6x+1) \\
& \\
&\omega_4(x)=x(24x^3+36x^2+14x+1),\qquad \omega_5(x)=x(120x^4+240x^3+150x^2+30x+1).
\end{align*}
In the case of $Q(n,y,p)$ we have to deal with minus signs and the different powers of $z$. The polynomials then are given by
\[
\sum_{m=1}^k (-1)^{1-m}S(k,m)m!z^{k-m}.
\]
The factors in front of the polynomials satisfy the Stirling numbers of the first kind. We get for $Q(n,y,p)$

\begin{align}
Q(n,y,p)&=\sum_{k=0}^\infty(-p\, k)_n(-y)^k \nonumber \\
&=(-1)^{n+1}\left(\dfrac{y}{y+1}\right)\sum_{k=1}^n s(n,k)(-p)^k
\sum_{m=1}^k S(k,m)\,m!\left(-\dfrac{1}{y+1}\right)^m \label{9} \\
&=(-1)^{n+1}\left(\dfrac{y}{y+1}\right)\sum_{k=1}^n s(n,k)(-p)^k
\omega_k\left(-\dfrac{1}{y+1}\right). \nonumber
\end{align}

To prove \eqref{9} we use a general transformation in \cite{7}. In this paper Boyadzhiev proves the following transformation formula \cite[(1.1)]{7}
\[
\sum_{n=0}^\infty\dfrac{g^{(n)}(0)}{n!}f(n)\,x^n=
\sum_{n=0}^\infty\dfrac{f^{(n)}(0)}{n!}\sum_{k=0}^\infty S(n,k)x^k\,g^{(k)}(x).
\]
He applies this formula for the function
\[
g(y)=\dfrac{1}{1-y} \qquad |y|<1
\]
which gives
\[
\sum_{k=0}^\infty f(k)\,y^k=
\dfrac{1}{1-y}\sum_{n=0}^\infty \dfrac{f^{(n)}(0)}{n!}\omega_n\left(\dfrac{y}{1-y}\right).
\]
Application to \eqref{9} gives
\[
\sum_{k=0}^\infty(-p\, k)_n(-y)^k=
\dfrac{1}{1+y}\sum_{k=0}^\infty D^{(k)}
\Big[(-p\, y)_n\Big]_{y=0}\dfrac{1}{k!}\omega_k\left(-\dfrac{y}{1+y}\right).
\]
For the derivative we get
\[
 D^{(k)}\Big[(-p\, y)_n\Big]_{y=0}=(-1)^n k!\,s(n,k)\,p^k \qquad 0 \leq k \leq n.
\]
Substitution gives
\begin{align*}
Q(n,y,p)&
=\sum_{k=0}^\infty(-p\, k)_n(-y)^k \\
&=(-1)^n\dfrac{1}{1+y}\sum_{k=0}^n
s(n,k)\,p^k\,\omega_k\left(-\dfrac{y}{1+y}\right)
\qquad\qquad\qquad\qquad\qquad\qquad\qquad
\end{align*}
\begin{align*}
&=(-1)^n\dfrac{1}{1+y}\sum_{k=0}^n
s(n,k)\,p^k\sum_{m=0}^k S(k,m)\,m!\left(-\dfrac{y}{y+1}\right)^m \\
&=(-1)^{n+1}\dfrac{y}{1+y}\sum_{k=0}^n
s(n,k)\,(-p)^k\sum_{m=0}^k S(k,m)\,m!\left(-\dfrac{1}{y+1}\right)^m.
\end{align*}
This proves \eqref{9}.

\

Using the definition of the Pochhammer symbol we get for $Q(n,y,p)$
\[
Q(n,y,p)=p^n\sum_{k=0}^\infty(-y)^k \prod_{m=0}^{n-1}\left(k-\dfrac{m}{p}\right).
\]
This summation of a product can be transformed into a hypergeometric function
\[
Q(n,y,p)=(-1)^{n-1}\dfrac{p\, \Gamma(n)}{y^{2n}}\ _n F_{n-1}\left(
\begin{array}{c}
	1-\dfrac{n-1}{p},\dots,1-\dfrac{1}{p},2 \\
	-\dfrac{n-1}{p},\dots,-\dfrac{1}{p}
\end{array}
;-y\right).
\] 
$_n F_{n-1}(x)$ is the generalized hypergeometric function.

\

The hypergeometric function can be written as a polylogarithm function of negative order. We get
\begin{equation}
Q(n,y,p)=(-1)^n\sum_{k=0}^n s(n,k)Li_{-k}(-y)\, p^k \quad n>0.
\label{10}
\end{equation}

\

It can be shown that $Q(n,y,p)$ satisfies a diferential-difference equation
\[
Q(n+1,y,p)=-p\, y\dfrac{d}{dy}Q(n,y,p)+n\, Q(n,y,p) \qquad \text{with}\qquad Q(0,y,p)=\dfrac{1}{(y+1)}.
 \]

At last $Q(n,y,p)$ can be written as a summation with an $n-$times differentiation of a power function:
 \[
Q(n,y,p)=\sum_{k=0}^\infty(-y)^{-p\, k+k+n}\dfrac{d^n}{dy^n}\left((-y)^{p\, k}\right). 
\]
Working out the derivative proves this property.

\

\

$\text{\bf\large Acknowledgments}$

\

I thank Prof. T.H. Koornwinder who pointed me to formula \cite[2.1.3 (12)]{8} for calculating the integral of the Catalan number of the second kind.

\

\end{document}